\documentclass[12pt]{amsart}
\usepackage{amssymb}

\title {On Parallel Lines and Free Group}

\author {Kwai-Man Fan}
\email {kmfan@math.ccu.edu.tw}
\date{May, 2007 }

\begin{document}

\maketitle

\abstract {Let $\Sigma=L_1 \cup \cdots \cup L_n \subset \Bbb C^2$
be a union of $n\ge 2$ distinct complex lines. Suppose $\pi_1(\Bbb
C^2\setminus \Sigma)$ is a free group. We show that $\Sigma$ is a
union of parallel lines.}
\endabstract

\section {Introduction}

Let $\Sigma=L_1 \cup \cdots \cup L_n \subset \Bbb C^2$ be a union
of $n\ge 2$ distinct lines. Here, a line means the solution set of
a linear equation of the form $ax+by+c=0$ in the complex plane
$\Bbb C^2$, where $a,b,c$ are complex numbers. In case $a,b,c$ are
real numbers, we may imagine a picture of the lines of $\Sigma$ in
the real plane $\Bbb R^2$. If these lines have the same slope,
$\Sigma$ is called a union of parallel lines. Note  $\Sigma
\subset \Bbb C^2$ is a union of parallel lines if and only if
$\Sigma$ has no singular point. Consider the four manifold $M=\Bbb
C^2\setminus \Sigma$. Suppose $\Sigma$ is a union of $n$ distinct
parallel lines. By performing a linear transformation, we may
assume that each of the line $L_i$ is given by an equation of the
form $x=c_i$, and we have $c_i\ne c_j$ if $i\ne j$. It follows
that $M \cong \Bbb C\setminus \{c_1\cdots,c_n\} \times \Bbb C$.
Hence we have $\pi_1(\Bbb C^2 \setminus \Sigma) \cong \pi_1(\Bbb
C\setminus \{c_1, \cdots, c_n\}) \cong F_n$, where $F_n$ is a free
group of rank $n$. For a union of lines $\Sigma$, in general, it
is difficult to decide the structure of the fundamental group
$\pi_1(\Bbb C^2 \setminus \Sigma)$ [8]. It is interesting to know
when the fundamental group $\pi_1(\Bbb C^2 \setminus \Sigma)$ will
decide the structure of $\Sigma$. In [2], the authors show that
$\pi_1(\Bbb C^2 \setminus \Sigma)$ is a free abelian group if and
only if $\Sigma$ is a union of lines in general position. Note a
union of parallel lines has one of the simplest configuration. One
question is: for a union of complex lines in the complex plane,
when is $\pi_1(\Bbb C^2 \setminus \Sigma)$ a free group? In this
paper, we prove:

\

\noindent {\bf Theorem 1.} Let $\Sigma=L_1 \cup \cdots \cup L_n
\subset \Bbb C^2$ be a union of $n\ge 2$ distinct complex lines.
Suppose $\pi_1(\Bbb C^2\setminus \Sigma)$ is a free group. Then
$\Sigma$ is a union of parallel lines.

\

\noindent Hence we have:

\

\noindent {\bf Corollary 1.} Let $\Sigma=L_1 \cup \cdots \cup L_n
\subset \Bbb C^2$ be a union of $n\ge 2$ distinct complex lines.
Then $\Sigma$ is a union of parallel lines if and only if
$\pi_1(\Bbb C^2\setminus \Sigma)$ is free.

\

\noindent This gives us a way to characterize a union of parallel
lines in the complex plane in terms of the fundamental group of
its complement. In section 3, we show:

\

\noindent {\bf Theorem 2.} Let $\Sigma=L_1 \cup \cdots \cup L_n
\subset \Bbb C^2$ be a union of $n\ge 2$ distinct complex lines.
Suppose $\Sigma$ has at least one singular point. Then $\pi_1(\Bbb
C^2\setminus \Sigma)$ contains a subgroup $H$ which is isomorphic
to a free abelian group of rank 2.

\

\noindent It is well known that a subgroup of a free group is free
[10], [11]. If $\Sigma$ is not a union of parallel lines, then
Theorem 2 implies that $\pi_1(\Bbb C^2\setminus \Sigma)$ is not a
free group. In order to prove Theorem 2, we need to study a
presentation of $G=\pi_1(\Bbb C^2\setminus \Sigma)$. Several
authors have given different ways to obtain a presentation of $G$,
c. f. for example, [1], [3], [9], [12], [14]. Here, we apply the
method given by Arvola in [1] to obtain a presentation of $G$. A
detail exposition of this method can be found, for example, in
Chapter 5 of the book [13] by Orlik and Terao. In this work, we
follow notations as given in [1], and in Section 5.3 of [13]. In
section 2, we give a summary of this method and derive a lemma.
The proof of Theorem 2 is given in section 3. Since the complex
projective plane is obtained from $\Bbb C^2$ by adjoining a
projective line, Corollary 1 implies:

\

\noindent {\bf Theorem 3.} Let $\Sigma=L_1 \cup \cdots \cup L_n
\subset \Bbb CP^2$ be a union of $n$ complex projective lines.
Then $\pi_1(\Bbb CP^2\setminus \Sigma)$ is a free group if and
only if $L_1 \cap \cdots \cap L_n \ne \emptyset$.

\

\noindent Theorem 3 and the Zariski's Hyperplane Section Theorem
[15] together implies:

\

\noindent {\bf Theorem 4.} Let $\Lambda=H_1 \cup \cdots \cup H_n
\subset \Bbb CP^{m+2}$ be a union of hyperplanes, where $m\ge 0$.
Then $\pi_1(\Bbb CP^{m+2}\setminus \Lambda)$ is a free group if
and only if $\dim(H_1  \cap \cdots \cap H_n)\ge m$.

\

\section {Some notes}

Let $\Sigma=L_1 \cup \cdots \cup L_n \subset \Bbb C^2$ be a union
of $n$ distinct complex lines. Let us begin with a summary of the
method of [1] which gives us a way to obtain a presentation of
$G=\pi_1(\Bbb C^2\setminus \Sigma)$. Choose coordinates $z_1,z_2$
for $\Bbb C^2$ and coordinates $x_1,y_1,x_2,y_2$ for $\Bbb R^4$.
Identify $\Bbb C^2$ with $\Bbb R^4$ by $z_1=x_1+i y_1$, and
$z_2=x_2+i y_2$. Let $\Bbb R^2$ be the span of $x_1,x_2$, and let
$\Bbb R^3$ be the span of $x_1,x_2,y_2$. Let $\phi^2: \Bbb R^4
\longrightarrow \Bbb R^2$, $\phi^3: \Bbb R^4 \longrightarrow \Bbb
\Bbb R^3$, and $\phi: \Bbb R^3 \longrightarrow \Bbb R^2$ be the
projection map given by $\phi^2(x_1,y_1,x_2,y_2)=(x_1,x_2)$,
$\phi^3(x_1,y_1,x_2,y_2)=(x_1,x_2,y_2)$, and
$\phi(x_1,x_2,y_2)=(x_1,x_2)$, respectively. Let $P$ be the set of
all multiple points of $\Sigma$. After a suitable linear
transformation, we assume the following two conditions holds:

\noindent Assumption $1)$ No lines of $\Sigma$ is of the form
$\{z_1=c\}$.

\noindent Assumption $2)$ If $p,p^{'} \in P$ are distinct multiple
points, then $x_1(p) \ne x_1(p^{'})$.

\noindent Write $P=\{p_1,\cdots,p_r\}$ so that
$x_1(p_1)<\cdots<x_1(p_r)$. Choose a piecewise linear  map \\$f :
\Bbb R \longrightarrow \Bbb R$ so that for each $p_i\in P$, the
graph of $f$ goes through $(x_1(p_i),y_1(p_i))$.   The map $f$ is
chosen to be a constant function at a small open neighborhood of
each $x_1(p_i)$. Such an $f$ is called a graphing map.

For a fixed $t$, let $K_t(f)=\{q \in \Bbb R^4\ |\ x_1(q)=t,
y_1(q)=f(t)\}$. Note $K_t(f)$ is a  complex line in $\Bbb C^2$
with the first coordinate $z_1=(t,f(t))$.  Let $\Gamma_f^4=\Sigma
\cap (\cup_{t\in \Bbb R} K_t(f))$. Conditions (1) and (2) given
above implies that for each $L_i\subset\Sigma$, $L_i\cap
\Gamma_f^4$ is homeomorphic to the real line $\Bbb R$, and
$\Gamma_f^4$ may be seen as a singular braid sitting in the
3-space $\cup_{t\in \Bbb R} K_t(f)\cong \Bbb R \times \Bbb R^2$.
Let $\Gamma_f^3=\phi^3(\Gamma_f^4)$, and
$\Gamma_f^2=\phi^2(\Gamma_f^4)$. By a suitable choice of
coordinates and $f$, we assume that the graph $\Gamma_f^2$ is
regular with respect to the projection map $\phi : (\Bbb
R^3,\Gamma_f^3) \longrightarrow (\Bbb R^2,\Gamma_f^2)$ in the
sense that $\phi$ and $f$ satisfy conditions of Lemma 5.62 of
[13]. For $\Gamma_f^2=\phi (\Gamma_f^3)$, the condition of $\phi$
being regular agrees with the usual definition as given in the
theory of singular braid. For $p\in P$, where $P$ is the set of
singular points of $\Sigma$, $\phi(P)$ is called an actual vertex
of $\Gamma_f^2$. For $q \in \Gamma_f^2$ such that
$|\phi^{-1}(q)\cap \Gamma_f^3|=2$, $q$ is called a virtual vertex
of $\Gamma_f^2$. Following the convention of Definition 5.68 of
[13], we assign a positive or negative sign to a virtual vertex of
$\Gamma_f^2$. The method of assigning a positive or negative sign
to a virtual vertex of $\Gamma_f^2$ agrees with the usual way that
is used in knot theory. Think that each of the string of the
singular braid $\Gamma_f^3 \subset \Bbb R^3$ is oriented from the
left side $(x_1=-\infty)$ to the right side $(x_1=+\infty)$
following the orientation of the $x_1$ axis. With respect to the
projection map $\phi : (\Bbb R^3,\Gamma_f^3) \longrightarrow (\Bbb
R^2,\Gamma_f^2)$, the positive virtual vertices are the positive
crossings, and the negative virtual vertices are the negative
crossings, c. f. Fig 5.11 of [13]. Call a regular 2-graph
$\Gamma_f^2$ admissible if its graphing map $f$  satisfies the
conditions of Lemma 5.66 of [13], and its virtual vertices are
marked by $+$ or $-$ signs accordingly. The graphing map $f$ is
chosen so that each of the vertical line $\{x_1=t\}$ in $\Bbb R^2$
goes through at most one vertex of $\Gamma_f^2$.

A presentation of $\pi_1(\Bbb C^2\setminus \Sigma)$ can be read
from an admissible 2-graph. Let $u^v=v^{-1} u v$, and let
$[u_1,\cdots,u_h]$ denote the cyclic relation:
$$u_1 u_2 \cdots u_h = u_2 \cdots u_h u_1=\cdots = u_h u_1 \cdots
u_{h-1}.$$ Let $\{q_1,\cdots,q_s\}$ be the set of all vertices of
$\Gamma_f^2$ so that $x_1(q_1)<\cdots< x_1(q_s)$. Choose $t_0<t_1<
\dots<t_{s-1}<t_s$ to separate the $x_1$ coordinates so that
$$t_0<x_1(q_1)<t_1<x_1(q_2)<t_2<\cdots<
x_1(q_{s-1})<t_{s-1}< x_1(q_s)<t_s.$$ Let $E(i) \subset \Bbb R^2$
be the vertical line given by $x_1=t_i$. Let $A(i)=E(i)\cap
\Gamma_f^2$, and we order the points of $A(i)$ by the $x_2$
coordinates in increasing order. Write $A(i)=(a_1^i,\cdots,a_n^i)$
so that $x_2(a_1^i)<x_2(a_2^i)<\cdots<x_2(a_n^i)$. For $i=0$, we
associate the word $g_j$ to the point $a_j^0$; for $j=1,\cdots,n$.
Let $W_0=(g_1,\cdots,g_n)$. The words $g_1,\cdots,g_n$ form a set
of generators for $G$. We associate a word
$w_j^i=w_j^i(g_1,\cdots,g_n)$ to the point $a_j^i$. Let
$W_i=(w_1^i,\cdots,w_n^i)$; for $i=1,\cdots,s$. For $i>0$, the set
of words $W_i$ are obtained from $W_{i-1}$ according to the rules
given below. Let $R(t_{i-1},t_i)=\{(x_1,x_2)\in \Bbb R^2\ |\
t_{i-1} \le x_1 \le t_i\}$. Let $q=q_i$ be the vertex of
$\Gamma_f^2$ in $R(t_{i-1},t_i)$, so that the segments of
$\Gamma_f^2$ in $R(t_{i-1},t_i)$ begin from
$a_j^{i-1},a_{j+1}^{i-1},\cdots,a_k^{i-1}$ intersect at the point
$q$, and the other segments of $\Gamma_f^2$ in $R(t_{i-1},t_i)$ do
not go through $q$. Suppose
$W_{i-1}=(w_1,\cdots,w_{j-1},w_j,\cdots,w_k,w_{k+1},\cdots,w_n)$.
Then $W_i=(w_1,\cdots,w_{j-1},w_k^{'},w_{k-1}^{'}
\cdots,w_{j+1}^{'},w_j^{'},w_{k+1},\cdots,w_n)$. There are three
cases:

\

\noindent Case 1. Suppose $q$ is an actual vertex. Then
$w_j^{'}=w_j,\ w_{j+1}^{'}=w_{j+1}^{w_j},\cdots,
w_k^{'}=w_k^{w_{k-1} \cdots w_j}$.

\

\noindent Case 2. Suppose $q$ is a positive virtual vertex. Then
$k=j+1$, $w_j^{'}=w_j$, $w_{j+1}^{'}=w_{j+1}^{w_j}$.

\

\noindent Case 3. Suppose $q$ is a negative virtual vertex. Then
$k=j+1$, $w_j^{'}=w_j^{w_{j+1}^{-1}}$, $w_{j+1}^{'}=w_{j+1}$.

\

\noindent To each of the actual vertex $q$ of $\Gamma_f^2$, we
associate a set of relations $R_q=[w_k,w_{k-1},\cdots,w_j]$.

\

\noindent {\bf Arvola's Theorem.} Let $\Sigma$ be an arrangement
of complex lines in $\Bbb C^2$. Let $M=\Bbb C^2\setminus \Sigma$.
Let $\Gamma$ be an admissible 2-graph for $\Sigma$. Then $\pi_1
(M) \cong\ <g_1,\cdots,g_n\ |\ \cup_q R_q \ >$ where $q$ ranges
over all the actual vertices of $\Gamma$.

\

\noindent If we range through a vertex $q_i$ of $\Gamma$, we need
to change the words in $W_i$ accordingly. Remark 4.10 of [1] is
helpful to us which says that we may begin to derive a
presentation of $\pi_1 (M)$ from the first actual vertex. Suppose
$q_1,\cdots,q_m$ are virtual vertices and $q_{m+1}$ is the first
actual vertex of an admissible 2-graph $\Gamma_f^2$ for $\Sigma$.
Let $W_m=(g_1(m),\cdots,g_n(m))$.  Suppose we range through a
positive virtual vertex. As in case 2, we get $w_j^{'}=w_j$ and
$w_{j+1}^{'}=w_j^{-1} w_{j+1} w_j$. Reverse these relations, we
get $w_{j+1}=w_j^{'} w_{j+1}^{'} {w_j^{'}}^{-1}$. Suppose we range
through a negative virtual vertex. As in case 3, we get
$w_j^{'}=w_{j+1} w_j w_{j+1}^{-1}$ and $w_{j+1}^{'}=w_{j+1}$.
Reverse these relations, we get $w_j={w_{j+1}^{'}}^{-1} w_j^{'}
w_{j+1}^{'}$. Hence we may express each $g_i$ of $W_0$ in terms of
$g_1(m),\cdots,g_n(m)$. By applying Tietze transformations [10],
we may replace $g_1,\cdots,g_n$ by $g_1(m),\cdots,g_n(m)$ as a set
of generators for $G$. In the next section, we are going to choose
generators for $G$ at $E(m)$.

Let $G \cong\ <g_1,\cdots,g_n\ |\ \cup_q R_q \ >$ be a
presentation of $\pi_1 (M)$ given as above. Suppose $q=q_i$ is an
actual vertex of $\Gamma$. Let
$W_i=(w_1,\cdots,w_j,\cdots,w_k,\cdots,w_n)$. Suppose that
$R_q=[w_k, w_{k-1},\cdots,w_j]$. $R_q$ is equivalent to the
following $k-j$ relations: $$w_k w_{k-1} \cdots w_j=w_{k-1} \cdots
w_j w_k,$$ $w_{k-1} w_{k-2}\cdots w_j w_k=w_{k-2} \cdots w_j w_k
w_{k-1}$, $\cdots$, $w_{j+1}w_j w_k \cdots w_{j+2}=w_j w_k \cdots
w_{j+2} w_{j+1}$. These $k-j$ relations are of form $ab=ba$.
Rewrite these relations in the form of $a b a^{-1} b^{-1}=e$,
where $e$ is the identity element of $\pi_1 (M)$. For each actual
vertex $q$ of $\Gamma$, we replace $R_q$ by relations of form
$a(q) b(q) a(q)^{-1} b(q)^{-1}=e$. Note $a(q)$ and $b(q)$ are
words in $g_1,\cdots,g_n$. Hence the above presentation of $\pi_1
(M)$ has the form:
$$\pi_1 (M) \cong\ <g_1,\cdots,g_n\ |\ U_1,\cdots,U_t>,$$ so that
each of the relator $U_i$ is a commutator [3]. That is: $U_i=a_i
b_i a_i^{-1} b_i^{-1}$, where $a_i,b_i$ are words in $g_1,\cdots
g_n$; for $i=1,\cdots,t$.

Let $H \cong <b_1,\cdots,b_n\ |\ V_1,\cdots,V_t>$ be a group so
that for $j=1,\cdots,t$, the relator $V_j$ is a commutator. Let
$W=W(b_1,\cdots,b_n)$ be a word in $b_1,\cdots,b_n$. We give a
condition so that $W$ represents a nontrivial element of $H$.
Recall the definition of exponent sum of a word $W$ as given in
pp. 76 of [10]. Suppose $W=b_{\upsilon_1}^{\alpha_1} \cdots
b_{\upsilon_r}^{\alpha_r}$, where the $\alpha_i$ are integers, and
$\upsilon_i\in \{1,\cdots,n\}$. Then the exponent sum of $W$ in
$b_\upsilon$ is the integer $\sigma_{\upsilon}=\displaystyle
\sum_{\upsilon_i=\upsilon} \alpha_i$. For the empty word $e$, let
$\sigma_{\upsilon}(e)=0$. For example, suppose $W=b_1^5 b_2
b_1^{-1} b_2^{-3} b_1^{-1}$. Then $\sigma_1(W)=3$,
$\sigma_2(W)=-2$. The exponent sum $\sigma_{\upsilon}$ has the
property, c. f. pp. 76 of [10], that for any word $W=U V$, we have
$$\sigma_{\upsilon}(W
)=\sigma_{\upsilon}(U) + \sigma_{\upsilon}(V). \eqno (1)$$

\

\noindent {\bf Lemma 1.} Let $H \cong <b_1,\cdots,b_n\ |\
V_1,\cdots,V_t>$ be a group, so that the relator $V_j$ is a
commutator; for $j=1,\cdots,t$. Let $A,B$ be two words in
$b_1,\cdots,b_n$. Suppose $A=B$ holds in $H$. Then we have
$\sigma_i(A)=\sigma_i(B)$; for $i=1,\cdots,n$. In particular,
suppose $A=e$, where $e$ is the identity element of $H$. Then we
have $\sigma_1(A)=\cdots=\sigma_n(A)=0$.

\

\noindent {\it Proof.} Recall that in a group $G$ with generators
$a,b,\cdots$, the trivial relators $a a^{-1}$, $a^{-1} a$, $b
b^{-1}$, $b^{-1} b, \cdots$ represent the identity element.
Suppose $A=B$ holds in $H$. Then we can change the word $A$ to $B$
by a finite number of steps $A=A_1,A_2\cdots,A_k=B$ so that
$A_{i+1}$ is obtained from $A_i$ by applying one of the following
two operations:

\noindent $i)$ Insertion of one of the words
$V_1,V_1^{-1}\cdots,V_t,V_t^{-1}$ or one of the trivial relators
of $H$ between any two consecutive symbols of $A_i$, or before
$A_i$, or after $A_i$.

\noindent $ii)$ Deletion of one of the words
$V_1,V_1^{-1}\cdots,V_t,V_t^{-1}$ or one of the trivial relators
of $H$ between any two consecutive symbols of $A_i$, or before
$A_i$, or after $A_i$.

\noindent Since $V_j$ is a commutator, it follows that
$\sigma_1(V_j)=\cdots=\sigma_n(V_j)=0$. For any trivial relator of
the form $g g^{-1}$, where $g\in
\{b_1,b_1^{-1}\cdots,b_n,b_n^{-1}\}$, we have $\sigma_1(g
g^{-1})=\cdots=\sigma_n(g g^{-1})=0$. By the property (1) of
exponent sum given above, we have
$$\sigma_i(A_1)=\sigma_i(A_2)=\cdots=\sigma_i(A_k).$$
Hence we have $\sigma_i(A)=\sigma_i(B)$; for $i=1,\cdots,n$.

Suppose  $A=e$, where $e$ is the identity element of $H$. Since
$\sigma_i(e)=0$, so we have $\sigma_i(A)=0$; for $i=1,\cdots,n$.
$\qed$

\

\noindent Lemma 1 has the following obvious implication. Let
$W=W(b_1,\cdots,b_n)$ be a word in $b_1,\cdots,b_n$. Suppose there
is an $i\in \{1,\cdots,n\}$ such that the exponent sum
$\sigma_i(W)\ne 0$. Then in the group $H$, we have $W \ne e$,
where $e$ is the identity element of $H$.

\

\section {Proof of main results}

Let us give the following:

\

\noindent {\it Proof of Theorem 2.} Choose suitable coordinates,
we assume $L_1,\cdots,L_n$ satisfy the two assumptions listed in
section 2. Choose a suitable graphing map $f : \Bbb R
\longrightarrow \Bbb R$ so that the graph $\Gamma=\Gamma_f^2$ is
regular with respect to the projection map $\phi : (\Bbb
R^3,\Gamma_f^3) \longrightarrow (\Bbb R^2,\Gamma_f^2)$. We assume
that the vertices and nodes of $\Gamma$ satisfy the conditions of
Lemma 5.66 of [13]. The graphing map $f$ is chosen so that each of
the vertical line $\{x_1=t\}$ in $\Bbb R^2$ goes through at most
one vertex of $\Gamma$. Following the convention of Definition
5.68 of [13], we assign a positive or negative sign to each
virtual vertex of $\Gamma$. Let $\{q_1,\cdots,q_s\}$ be the set of
all vertices of $\Gamma$ so that $x_1(q_1)<\cdots< x_1(q_s)$.
Choose $t_0<t_1< \dots<t_{s-1}<t_s$ to separate the $x_1$
coordinates of these vertices so that
$$t_0<x_1(q_1)<t_1<x_1(q_2)<t_2<\cdots<
x_1(q_{s-1})<t_{s-1}< x_1(q_s)<t_s.$$ Let $E(i)$ be the vertical
line $\{x_1=t_i\} \subset \Bbb R^2$. Let $q_{m+1}$ be the first
actual vertex of $\Gamma$. Let $E(m)\cap
\Gamma=\{(t_m,y_1),\cdots,(t_m,y_n)\}$, so that we have
$y_1<\cdots<y_n$. Let $Q_i=(t_m,y_i)$, and we associate the symbol
$g_i$ to the point $Q_i$; for $i=1,\cdots,n$. Following Remark
4.10 of [1], we may use $g_1,\cdots,g_n$ as a set of generators
for $G$. Let
$$\pi_1 (M) \cong\ <g_1,\cdots,g_n\ |\ \cup_q R_q \ >$$
be the presentation of $G$ given by Theorem 4.7 of [1], using
$g_1,\cdots,g_n$ as generators. Suppose the first actual vertex
$q=q_{m+1}$ of $\Gamma$ gives us the relation
$R_q=[g_k,g_{k-1},\cdots,g_j]$ so that $k-j\ge 1$. Write this
presentation, as given in section 2, in the form:
$$G=\pi_1 (M) \cong\ <g_1,\cdots,g_n\ |\ U_1,\cdots,U_t>$$
so that $U_i$ is a commutator; for $i=1,\cdots,t$. In this
presentation of $\pi_1 (M)$, the relation $R_q$ gives us $k-j$
relations of form $aba^{-1}b^{-1}=e$, where $e$ is the identity
element of $G$. The first one of these $k-j$ relations gives us:
$$g_k (g_{k-1}\cdots g_j)=(g_{k-1}\cdots g_j) g_k. \eqno (2)$$
Let $A=g_k$ and $B=g_{k-1}\cdots g_j$. Note (2) implies that the
relation $AB=BA$ holds in $G$. Let $H$ be the subgroup of $G$
generated by $A$ and $B$. Hence we have $H=\{A^i B^j \ |\ i,j \in
\Bbb Z\}$. Suppose $A^{i_0}B^{j_0}=A^{i_1}B^{j_1}$. Then we have
$$W=A^{i_0}B^{j_0} (A^{i_1}B^{j_1})^{-1} = A^{i_0-i_1}B^{j_0-j_1}
=g_k^{i_0-i_1}(g_{k-1}\cdots g_j)^{j_0-j_1}=e.$$
Consider the exponent sum
$\sigma_j(W),\cdots,\sigma_{k-1}(W),\sigma_k(W)$. Note
$\sigma_k(W)=i_0-i_1$, and
$$\sigma_j(W)=\cdots=\sigma_{k-1}(W)=j_0-j_1.$$
Since $W=e$ holds in $G$, Lemma 1 implies that $i_0-i_1=0$ and
$j_0-j_1=0$. Hence we have $A^{i_0}B^{j_0}=A^{i_1}B^{j_1}$ if and
only if $i_0=i_1$ and $j_0=j_1$ holds. In particular, this implies
that for the subgroups $<A>, <B>$ of $G$ generated by $A,B$
respectively, we have $<A> \cap <B> =\{e\}$ and $<A> \cong <B>
\cong \Bbb Z$. Hence we have $H \cong \Bbb Z \oplus \Bbb Z$.
$\qed$

\

Let us give the following:

\

\noindent {\it Proof of Theorem 1.} Let $\Sigma=L_1 \cup \cdots
\cup L_n \subset \Bbb C^2$ be a union of $n\ge 2$ distinct complex
lines. Suppose $\Sigma$ is not a union of parallel lines. Then
$\Sigma$ has least one singular point. Theorem 2 implies that
$G=\pi_1(\Bbb C^2\setminus \Sigma)$ contains a subgroup $H$ which
is isomorphic to a free abelian group of rank 2. Since a subgroup
of a free group is free, it follows that $G$ is not a free group.
This shows that if $\pi_1(\Bbb C^2\setminus \Sigma)$ is a free
group, then $\Sigma \subset C^2$ is a union of parallel lines.
$\qed$

\

\noindent {\it Proof of Corollary 1.} If $\Sigma= L_1 \cup \cdots
\cup L_n \subset \Bbb C^2$ be a union of $n\ge 2$ parallel lines,
then we have $\pi_1(\Bbb C^2\setminus \Sigma) \cong F_n$. This
fact together with Theorem 1 gives us Corollary 1. $\qed$

\

\noindent {\it Proof of Theorem 3.} Let $\Sigma=L_1 \cup \cdots
\cup L_n \subset \Bbb CP^2$ be a union of $n$ distinct complex
projective lines. Suppose $n=1$. Then $\Bbb CP^2 \setminus L_1
\cong \Bbb C^2$ is simply connected. Assume $n\ge2$. Note for any
$L_i \subset \Sigma$, $\Bbb CP^2 \setminus L_i \cong \Bbb C^2$,
and $\Sigma \setminus L_i \subset \Bbb CP^2 \setminus L_i$ is a
union of $n-1$ distinct affine lines. Suppose $n=2$. Then $L_1
\cap L_2$ contains exactly one point, and $\Bbb CP^2 \setminus
\Sigma \cong (\Bbb C\setminus \{0\}) \times \Bbb C$. It follows
that $\pi_1(\Bbb CP^2\setminus \Sigma) \cong \Bbb Z$.

Suppose $n > 2$. Note we have $\pi_1(\Bbb CP^2\setminus
\Sigma)=\pi_1((\Bbb CP^2 \setminus L_1)\setminus (\Sigma\setminus
L_1))$. $\Sigma\setminus L_1$ is a union of $n-1$ parallel lines
in $\Bbb CP^2 \setminus L_1 \cong \Bbb C^2$  if and only if $L_1
\cap \cdots \cap L_n \ne \emptyset$ holds. By Corollary 1, $\Sigma
\setminus L_1 \subset \Bbb CP^2 \setminus L_1 \cong \Bbb C^2$ is a
union of $n-1$ parallel lines if and only if $\pi_1(\Bbb
CP^2\setminus \Sigma) \cong \pi_1((\Bbb CP^2 \setminus
L_1)\setminus (\Sigma\setminus L_1))\cong F_{n-1}$. This gives us
Theorem 3. $\qed$

\

\noindent {\it Proof of Theorem 4.} Let $\Lambda=H_1 \cup \cdots
\cup H_n \subset \Bbb CP^{m+2}$ be a union of $n$ distinct
hyperplanes. The case when $m=0$ is settled in Theorem 3, so we
assume $m\ge 1$. Suppose $n=1$. Then $\Bbb CP^{m+2} \setminus H_1
\cong \Bbb C^{m+2}$. Since $\Bbb C^{m+2}$ is simply connected, the
conclusion follows. Suppose $n=2$. Note we have $\dim(H_1\cap
H_2)=m$ and $\Bbb CP^{m+2} \setminus (H_1 \cup H_2) \cong (\Bbb
C\setminus \{0\}) \times \Bbb C^{m+1}$. It follows that
$\pi_1(\Bbb CP^{m+2}\setminus \Lambda) \cong \Bbb Z$. Suppose $n >
2$. For a set of hyperplanes $H_1,\cdots,H_n$ in the projective
space $\Bbb CP^{m+2}$, $H_1 \cap \cdots \cap H_n$ is either an
empty set, or a subspace of $\Bbb CP^{m+2}$ of dimension $k \le
m$.

\

\noindent $i)$ Suppose $H_1 \cap \cdots \cap H_n = E$, where $E
\cong \Bbb CP^m$ is a subspace of $\Bbb CP^{m+2}$ of dimension
$m$. Let $\Bbb P^2 \subset \Bbb CP^{m+2}$ be a 2-plane so that
$a)$ $\Bbb P^2$ is not a subspace of $H_i$; for $i=1,\cdots,n$,
and $b)$ $\Bbb P^2$ intersects $E$ transversally. It follows that
$\Bbb P^2 \cap H_i =L_i$ is a complex projective lines in $\Bbb
P^2$. Since $\dim(\Bbb P^2)+\dim(E)=2+m$, so $\Bbb P^2 \cap E$ has
exactly one point. Let $\Sigma=\Bbb P^2 \cap \Lambda = L_1 \cup
\cdots \cup L_n$. By the Zariski Hyperplane Section Theorem, we
have $\pi_1(\Bbb CP^{m+2}\setminus \Lambda) \cong \pi_1(\Bbb P^2
\setminus \Sigma)$.  Note $$L_1 \cap \cdots \cap L_n =(\Bbb P^2
\cap H_1) \cap \cdots \cap (\Bbb P^2 \cap H_n) =\Bbb P^2 \cap (H_1
\cap \cdots   \cap H_n) =\Bbb P^2 \cap E \ne \emptyset .$$
\noindent By Theorem 3, it follows that $\pi_1(\Bbb
CP^{m+2}\setminus \Lambda) \cong \pi_1(\Bbb P^2 \setminus  \Sigma)
\cong F_n$.

\

\noindent $ii)$ Suppose we have either $H_1 \cap \cdots \cap H_n=
\emptyset$, or  $\dim(H_1 \cap \cdots \cap H_n) < m$. Let $\Bbb
P^2 \subset \Bbb CP^{m+2}$ be a 2-plane so that the following
three conditions hold:

\noindent 1) For $1 \le i \le n$, $\dim(\Bbb P^2 \cap H_i)=1$.

\noindent 2) For $1 \le i < j \le n$, $\Bbb P^2 \cap H_i \cap H_j$
is a point.

\noindent 3) For $1 \le i < j <k \le n$, suppose $\dim(H_i \cap
H_j \cap H_k) < m$. Then $\Bbb P^2 \cap H_i \cap H_j \cap H_k =
\emptyset$.

\noindent Suppose $H_{i_1} \cap \cdots \cap H_{i_t} \ne \emptyset$
and $\dim ( H_{i_1} \cap \cdots \cap H_{i_t} ) < m$, where
$\{i_1,\cdots,i_t\}\subset \{1,\cdots,n\}$. Then there are $i,j,k
\in \{i_1,\cdots,i_t\}$ such that $\dim(H_i \cap H_j \cap H_k) <
m$. Suppose for any $i,j,k \in \{i_1,\cdots,i_t\}$, we have
$\dim(H_i \cap H_j \cap H_k) \ge m$. This implies $\dim ( H_{i_1}
\cap \cdots \cap H_{i_t} ) \ge m$ which is a contradiction to the
hypothesis.

Let $\Bbb P^2 \cap H_i =L_i$, and let $\Sigma=\Bbb P^2 \cap
\Lambda = L_1 \cup \cdots \cup L_n.$ By the Zariski Hyperplane
Section Theorem, we have $\pi_1(\Bbb CP^{m+2}\setminus \Lambda)
\cong \pi_1(\Bbb P^2 \setminus \Sigma)$. Note
$$L_1 \cap \cdots \cap L_n =(\Bbb P^2
\cap H_1) \cap \cdots \cap (\Bbb P^2 \cap H_n) =\Bbb P^2 \cap (H_1
\cap \cdots \cap H_n).$$ Since we have either $H_1 \cap \cdots
\cap H_n= \emptyset$ or $\dim(H_1 \cap \cdots \cap H_n) < m$, the
third condition implies that $L_1 \cap \cdots \cap L_n=\Bbb P^2
\cap (H_1 \cap \cdots \cap H_n) =\emptyset$. Theorem 3 implies
$\pi_1(\Bbb P^2 \setminus \Sigma) \ncong F_n$. It follows that
$\pi_1(\Bbb CP^{m+2}\setminus \Lambda) \ncong F_n$. $\qed$

\

\enddocument